\documentclass [leqno, 10pt]{amsart}

\usepackage[all]{xy}
\usepackage{amssymb, amsmath, amsfonts, amsthm, amsbsy, amscd, stmaryrd}
\usepackage[bbgreekl]{mathbbol} 
\usepackage[mathscr]{euscript}
\usepackage{esint}
\usepackage{upgreek}
\usepackage{enumerate}

\usepackage{url}
\usepackage[linktocpage]{hyperref}

\newtheorem*{theorem}{Theorem}
\newtheorem*{theorem*}{Theorem}{\,}

\theoremstyle{definition}
\newtheorem*{example}{Example}

\makeatletter
\def\l@subsection{\@tocline{2}{0pt}{1pc}{4.6em}{}}
\renewcommand{\tocsubsection}[3]{%
  \indentlabel{\@ifnotempty{#2}{\hspace*{2.3em}\makebox[2.3em][l]{%
    \ignorespaces#1 #2.\hfill}}}#3}
\makeatother

\def\cprime{$'$}

\newcommand{\Hol}{\text{Hol}}

\newcommand{\imt}{\iota}

\newcommand{\dsc}{\Lambda}

\newcommand{\symplecto}{\mathbb{symp}}

\newcommand{\Om}{\Omega}

\newcommand{\sflat}{\curlyvee}

\newcommand{\vol}{\mathsf{vol}}

\newcommand{\dum}{\,\cdot\,\,}
\newcommand{\Ga}{\Gamma}

\newcommand{\pr}{\partial}
\newcommand{\ctm}{T^{\ast}M}

\newcommand{\lie}{\mathfrak{L}}

\newcommand{\be}{\beta}

\newcommand{\rea}{\mathbb R}

\newcommand{\tr}{\operatorname{\mathsf{tr}}}

\begin{document}
\title{Infinitesimal affine automorphisms of symplectic connections}

\author{Daniel J. F. Fox} 
\address{Departamento de Matem\'aticas del \'Area Industrial \\ Escuela T\'ecnica Superior de Ingenier\'ia y Dise\~no Industrial\\ Universidad Polit\'ecnica de Madrid\\Ronda de Valencia 3\\ 28012 Madrid Espa\~na}
\email{daniel.fox@upm.es}

\begin{abstract}
Conditions are given under which an infinitesimal automorphism of a torsion-free connection preserving a symplectic form is necessarily a symplectic vector field. An example is given of a compact symplectic nilmanifold admitting a flat symplectic connection and an infinitesimal automorphism that is not symplectic.
\end{abstract}

\maketitle

On a symplectic manifold $(M, \Om)$, a \textit{symplectic connection} is a torsion-free affine connection $\nabla$ that preserves the symplectic form $\Om$. This note addresses the question of when an infinitesimal affine automorphism of a symplectic connection must be an infinitesimal symplectomorphism. Recall that an infinitesimal affine automorphism of an affine connection $\nabla$ is a vector field $X$ such that the Lie derivative $\lie_{X}\nabla$ vanishes. 

The automorphism group of the standard flat affine connection $\nabla$ on $n$-dimensional affine space is the full group of affine transformations. On the other hand, the automorphism group of a $\nabla$-parallel metric or symplectic form is the proper subgroup of isometric affine transformations or symplectic affine transformations. In particular, the full group of affine transformations preserves the standard flat connection, but only its subgroup comprising symplectic affine transformation preserves a parallel Darboux symplectic form. Similar examples on solvable symplectic Lie groups abound. This shows that some condition is necessary to conclude that an infinitesimal automorphism of a connection preserving some tensor necessarily preserves the tensor.

On the other hand, for the Levi-Civita connection $D$ of a Riemannian metric $g$ on a compact manifold, Theorem $4$ of K. Yano's \cite{Yano} shows that an infinitesimal automorphism of $D$ is a $g$-Killing field, and consequently that if a connected Lie group acts on a compact Riemannian manifold by automorphisms of the Levi-Civita connection, it acts by isometries. Related more general results are recorded in chapter VI of \cite{Kobayashi-Nomizu}. Yano's theorem means that, if $M$ is compact, the quotient of the affine automorphism group of $D$ by the isometry group of $g$ is discrete. For pseudo-Riemannian manifolds the analogous claims are false. In \cite{Boubel-Mounoud}, the compact Lorentzian $3$-manifolds admitting an affine automorphism that is not an isometry are classified. Many of the examples admit a one-parameter group of affine automorphisms that are not isometric. 

One expects the symplectic situation to more closely resemble the Lorentzian case than the Riemannian setting because the existence of isotropic subspaces means that holonomy need not act completely irreducibly. Here there are obtained results for infinitesimal affine automorphisms of symplectic connections analogous to those available for metric connections:

\begin{theorem}\label{liesymplectotheorem}
Let $(M, \Om)$ be a $2n$-dimensional symplectic manifold and let $\nabla$ be a torsion-free affine connection preserving $\Om$. Suppose the compactly supported vector field $X$ is an infinitesimal automorphism of $\nabla$. Then $X$ is symplectic, meaning $\lie_{X}\Om = 0$, in any of the following situations:
\begin{enumerate}
\item $M$ is noncompact. 
\item $\dim M =2$.
\item\label{euler} $M$ is compact with nonzero Euler characteristic.
\item $M$ is compact and the linear holonomy of $\nabla$ is irreducible.
\item $M$ is compact and $\nabla$ is the Levi-Civita connection of a K\"ahler metric compatible with a constant multiple of $\Om$.
\end{enumerate}
\end{theorem}

It will also be shown that the condition on the Euler characteristic in \eqref{euler} of the theorem cannot be removed. Since a compact nilmanifold is parallelizable it has vanishing Euler characteristic. Since there are many symplectic nilmanifolds (that, moreover, need not admit K\"ahler structures) it is natural to look among them for an example.
After the proof of the theorem, it is shown that there is a compact four-dimensional symplectic nilmanifold that admits a family of flat symplectic connections that each admit an infinitesimal automorphism that is not symplectic.

\begin{proof}[Proof of Theorem]
For a vector field $X \in \Ga(TM)$, let $X^{\sflat} \in \Ga(\ctm)$ be the symplectically dual one-form defined by $X^{\sflat} = \Om(X, \dum)$. 
Let $\symplecto(M, \Om) \subset \Ga(TM)$ be Lie subalgebra comprising symplectic vector fields, namely those vector fields preserving $\Om$. Since $\lie_{X}\Om = dX^{\sflat}$, $\symplecto(M, \Om)$ consists exactly in those vector fields whose symplectically dual one-form is closed. 

The abstract index conventions are used, and indices are raised and lowered (preserving horizontal position) using the symplectic form $\Om_{ij}$ and the dual antisymmetric bivector $\Om^{ij}$, consistently with the conventions $X_{i} = X^{p}\Om_{pi}$ and $X^{i} = \Om^{ip}X_{p}$ (so that $\Om^{ip}\Om_{pj} = -\delta_{j}\,^{i}$, where $\delta_{j}\,^{i}$ is the canonical pairing between the tangent space and its dual). Enclosure of indices in square brackets indicates complete antisymmetrization over the enclosed indices. A label is in either \textit{up} position or \textit{down}, and a label appearing as both an up and a down index indicates the trace pairing (summation convention). With these conventions, $X^{\sflat}_{i}$ and $X_{i}$ are synonyms, and $X^{p}Y_{p} =  -X_{q}Y^{q} = -X^{p}Y^{q}\Om_{pq} = \Om(X, Y)$. 

The curvature $R_{ijk}\,^{l}$ of a torsion-free affine connection $\nabla$ is defined by $2\nabla_{[i}\nabla_{j]}X^{k} = R_{ijp}\,^{k}X^{p}$ for $X \in \Ga(TM)$. 
By the Ricci identity, $0 = 2\nabla_{[i}\nabla_{j]}\Om_{kl} = -2R_{ij[kl]}$, in which $R_{ijkl} = R_{ijk}\,^{p}\Om_{pl}$. 

The Lie derivative along $X \in \Ga(TM)$ of a torsion-free affine connection $\nabla$ is given by
\begin{align}\label{liex}
(\lie_{X}\nabla)(Y, Z) = [X, \nabla_{Y}Z] - \nabla_{[X, Y]}Z - \nabla_{Y}[X, Z] = \nabla_{Y}\nabla_{Z}X - \nabla_{\nabla_{Y}Z}X + R(X, Y)Z.
\end{align}
Equivalently $(\lie_{X}\nabla)_{ij}\,^{k} = \nabla_{i}\nabla_{j}X^{k} + X^{p}R_{pij}\,^{k}$. If $\nabla$ is symplectic, then 
\begin{align}\label{lienabla0}
\begin{split}
(\lie_{X}\nabla)_{ijk} & = \nabla_{i}\nabla_{j}X_{k} + X^{p}R_{pijk}.
\end{split}
\end{align}
For any symplectic connection $\nabla$, there holds $dX^{\sflat}_{ij} = 2\nabla_{[i}X^{\sflat}_{j]}$.
Skew-symmetrizing \eqref{lienabla0} in $jk$ yields
\begin{align}
\begin{split}
0 & = 2(\lie_{X}\nabla)_{i[jk]} = 2\nabla_{i}\nabla_{[j}X_{k]} + 2X^{p}R_{pi[jk]} = \nabla_{i}dX^{\sflat}_{jk}.
\end{split}
\end{align}
so that $dX^{\sflat}$ is parallel. If $M$ is noncompact and $X$ is compactly supported, then $dX^{\sflat}$ vanishes outside the support of $X$, and so vanishes everywhere because it is parallel.

Contracting \eqref{lienabla0} yields $0 = \nabla_{i}\nabla_{p}X^{p}$, so that $\nabla_{p}X^{p}$ is constant. Write $\Om_{k} = \tfrac{1}{k!}\Om^{k}$.
If $M$ is compact, then 
\begin{align}
\begin{split}
0 &= \int_{M}d(X^{\sflat} \wedge \Om_{n-1}) = \int_{M}dX^{\sflat}\wedge \Om_{n-1} = \int_{M}(\nabla_{p}X^{p})\Om_{n} = (\nabla_{p}X^{p})\vol_{\Om_{n}}(M)
\end{split}
\end{align}
implies $\nabla_{p}X^{p} = 0$. Whether or not $M$ is compact, $dX^{\sflat} \wedge \Om_{n-1} = (\nabla_{p}X^{p})\Om_{n} = 0$. If $2n = 2$, this means $dX^{\sflat} = 0$, so that $X \in \symplecto(M, \Om)$.

Suppose $M$ is compact. Let $A_{i}\,^{j} = dX^{\sflat}_{i}\,^{j} = \nabla_{i}X^{j} - \nabla^{j}X_{i}$. For $k \geq 2$ define 
\begin{align}
A^{\circ k}\,_{i}\,^{j} = A_{i}\,^{i_{1}}A_{i_{1}}\,^{i_{2}}\dots A_{i_{k-2}}\,^{i_{k-1}}A_{i_{k-1}}\,^{j}, 
\end{align}
so that $A^{\circ k}\,_{i}\,^{j}$ is the $k$-fold composition of the endomorphism $A_{i}\,^{j}$. Because $A^{\circ k}\,_{ji} = - A^{\circ k}\,_{ij}$, the trace $\tr A^{\circ k}$ satisfies 
\begin{align}\label{trk}
\tr A^{\circ k} = A^{\circ k}\,_{p}\,^{p} = - A^{ij}A^{\circ k-1}\,_{ij} =  - A^{\circ k-1}\,_{ij} \nabla^{i}X^{j}.  
\end{align}
Since $A_{i}\,^{j}$ is parallel, $\tr A^{\circ k}$ is equal to some constant $c_{k} \in \rea$. Because $dX^{\sflat}_{ij} = A_{ij}$ is parallel, for $k \geq 2$, using \eqref{trk} and integrating by parts shows 
\begin{align}
c_{k}\vol_{\Om_{n}}(M) = \int_{M}\tr A^{\circ k}\,\Om_{n} = -2\int_{M}A^{\circ k-1}\,_{ij} \nabla^{i}X^{j}\,\Om_{n} = 2\int_{M}X^{j}\nabla^{i}A^{\circ k-1}\,_{ij} \,\Om_{n} = 0,
\end{align}
so $0 = c_{k} = \tr A^{\circ k}$. Since this holds for all $k \geq 2$, and $A_{p}\,^{p} = 0$, the endomorphism $A_{i}\,^{j}$ is nilpotent.

Now it will be shown that if $dX^{\sflat}$ is not identically zero then $\chi(M) = 0$. Since $A_{i}\,^{j}$ is parallel and nilpotent, the kernel and image of any power $A^{\circ k}$ are parallel nontrivial proper subbundles of $TM$. If $k$ is the largest integer such that $A^{\circ k} \neq 0$, then, since $A_{ij} = -A_{ji}$, the image $L$ of $A^{\circ k-1}$ is a parallel isotropic subbundle of $TM$. This is enough to conclude $\chi(M) = 0$, as follows. Let $J$ be an almost complex structure compatible with $\Om$, with corresponding Riemannian metric $g_{ij} = \Om_{ip}J_{j}\,^{p}$. Then $TM = L \oplus JL \oplus K$, where $K$ is the $g$-orthogonal complement of $L \oplus JL$. By Theorem $1.4$ of \cite{Dazord} there vanishes the Euler class of a symplectic vector bundle admitting a Lagrangian subbundle. The bundle $L \oplus JL$ is a symplectic vector bundle (with the symplectic structure the restriction of $\Om$), and $L$ is a Lagrangian subbundle of $L \oplus JL$, so the Euler class of $L \oplus JL$ vanishes. The Euler class of $TM$ equals the product of the Euler class of $K$ and the Euler class of $L \oplus JL$, so the Euler class of $TM$ vanishes.
Hence if $\chi(M) \neq 0$, $X \in \symplecto(M, \Om)$.

Let $\Hol_{x}$ be the linear holonomy group of $\nabla$ acting on $T_{x}M$. Since $A_{i}\,^{j}$ is parallel, as an endomorphism of $T_{x}M$ it commutes with $\Hol_{x}$. 
This implies that $\Hol_{x}$ preserves the filtration of $T_{x}M$ determined by the generalized null spaces of $A_{i}\,^{j}$. In particular, if the action of $\Hol_{x}$ on $T_{x}M$ is irreducible, then the kernel of $A_{i}\,^{j}$ must be all of $T_{x}M$, so that $A_{i}\,^{j} = 0$. This shows $X \in \symplecto(M,\Om)$ (this argument is the proof of $(1)$ of Theorem $1$ of Appendix $5$ of \cite{Kobayashi-Nomizu}).

If $M$ is compact and $\nabla$ is the Levi-Civita connection of a K\"ahler metric $g$ compatible with a constant multiple of $\Om$, then an infinitesimal automorphism $X$ of $\nabla$ is a $g$-Killing field by the already cited theorem of Yano, and $\Om$ is a harmonic form, so preserved by $X$.
\end{proof}

The theorem implies that if a symplectic connection on a compact symplectic manifold admits an infinitesimal automorphism that is not symplectic, then the manifold must have dimension at least $4$ and Euler characteristic zero and the holonomy of the symplectic connection must preserve a nontrivial filtration of the tangent bundle. 

\begin{example}
Let $M$ be the quotient of $\rea^{4}$ by the action of the discrete affine group $\dsc$ generated by unit translations along the $t^{1}$, $t^{2}$, and $t^{3}$ axes, and the map $(t^{1}, t^{2}, t^{3}, t^{4}) \to (t^{1} + t^{2}, t^{2}, t^{3}, t^{4} +1)$. The one-forms $e^{1} = dt^{1} - t^{4}dt^{2}$, $e^{2} = dt^{2}$, $e^{3} = dt^{3}$, and $e^{4} = dt^{4}$ are invariant under the action of $\dsc$. Since the symplectic form $\Om = e^{1} \wedge e^{2} + e^{3} \wedge e^{4}$ is $\dsc$-invariant, it descends to a symplectic structure on $M$. The symplectic manifold $(M, \Om)$ is one version of the Kodaira-Thurston nilmanifold. It is presented here as described in Example $2.1$ in chapter $2$ of \cite{Tralle-Oprea}. It is a $2$-torus bundle over a $2$-torus, so has Euler characteristic zero, and admits no K\"ahler metric.

The only nontrivial Lie bracket among the dual $\dsc$-invariant vector fields $E_{1} = \tfrac{\pr}{\pr t^{1}}$, $E_{2} = \tfrac{\pr}{\pr t^{2}} + t^{4}\tfrac{\pr}{\pr t^{1}}$, $E_{3} = \tfrac{\pr}{\pr t^{3}}$, and $E_{4} = \tfrac{\pr}{\pr t^{4}}$ is $[E_{2}, E_{4}] = -E_{1}$.
For any constant $\be\in \rea$, define the affine connection $\nabla$ by declaring 
\begin{align}
&\nabla_{E_{4}}E_{2} = -(\be - 2/3)E_{1},& &\nabla_{E_{2}}E_{4} = -(\be + 1/3)E_{1},& &\nabla_{E_{2}}E_{2} = -(\be+1/3)E_{3}, 
\end{align}
and all other covariant derivatives in the frame $E_{i}$ to be null. For $a \in \rea^{4}$ let $X^{a} = a^{1}E_{1} + a^{2}E_{3} + a^{3}E_{3} + a^{4}E_{4}$. Using
\begin{align}\label{christo}
\begin{split}
\nabla_{X^{a}}X^{b} &= \left((2/3 - \be)a^{4}b^{2} - (\be + 1/3)a^{2}b^{4}\right)E_{1} - (\be + 1/3)a_{2}b_{2}E_{3}\\
&= \left((2/3)a^{4}b^{2} - (1/3)a^{2}b^{4} - \be(a^{2}b^{4} + a^{4}b^{2})\right)E_{1} - (\be + 1/3)a_{2}b_{2}E_{3},
\end{split}
\end{align}
it is straightforward to check that $\nabla$ is torsion-free, symplectic, and flat. For example,
\begin{align}
\nabla_{X^{a}}X^{b} - \nabla_{X^{b}}X^{a} = -(a^{2}b^{4} - a^{4}b^{2})E_{1} = [X^{a}, X^{b}],
\end{align}
showing that $\nabla$ is torsion-free. That $\nabla \Om = 0$ follows because, by \eqref{christo}, $\Om(\nabla_{X^{a}}X^{b}, X^{c})$ is symmetric in $b$ and $c$. Since, by \eqref{christo}, $\nabla_{X^{a}}X^{b}$ is contained in the span of $E_{1}$ and $E_{3}$, and there vanish all covariant derivatives of the form $\nabla_{E_{i}}E_{j}$ in which either $i$ or $j$ equals either $1$ or $3$, there hold $\nabla_{X^{a}}\nabla_{X^{b}}X^{c} = 0$ and $\nabla_{\nabla_{X^{a}}X^{b}}X^{c} = 0$. Hence $\nabla$ is flat, and, by \eqref{liex}, any $X^{a}$ is an infinitesimal automorphism of $\nabla$.  Since $\lie_{X^{a}}\Om = d(\imt(X^{a})\Om) = -a^{2}e^{2} \wedge e^{4}$, the vector field $E_{2}$ is an infinitesimal automorphism of $\nabla$ that is not symplectic.
\end{example}

The results obtained here suggest that it is an interesting problem to classify symplectic connections admitting an affine automorphism that is not a symplectomorphism.

\bibliographystyle{amsplain}

\def\polhk#1{\setbox0=\hbox{#1}{\ooalign{\hidewidth
  \lower1.5ex\hbox{`}\hidewidth\crcr\unhbox0}}} \def\cprime{$'$}
  \def\cprime{$'$} \def\cprime{$'$}
  \def\polhk#1{\setbox0=\hbox{#1}{\ooalign{\hidewidth
  \lower1.5ex\hbox{`}\hidewidth\crcr\unhbox0}}} \def\cprime{$'$}
  \def\cprime{$'$} \def\cprime{$'$} \def\cprime{$'$} \def\cprime{$'$}
  \def\cprime{$'$} \def\polhk#1{\setbox0=\hbox{#1}{\ooalign{\hidewidth
  \lower1.5ex\hbox{`}\hidewidth\crcr\unhbox0}}} \def\cprime{$'$}
  \def\Dbar{\leavevmode\lower.6ex\hbox to 0pt{\hskip-.23ex \accent"16\hss}D}
  \def\cprime{$'$} \def\cprime{$'$} \def\cprime{$'$} \def\cprime{$'$}
  \def\cprime{$'$} \def\cprime{$'$} \def\cprime{$'$} \def\cprime{$'$}
  \def\cprime{$'$} \def\cprime{$'$} \def\cprime{$'$} \def\dbar{\leavevmode\hbox
  to 0pt{\hskip.2ex \accent"16\hss}d} \def\cprime{$'$} \def\cprime{$'$}
  \def\cprime{$'$} \def\cprime{$'$} \def\cprime{$'$} \def\cprime{$'$}
  \def\cprime{$'$} \def\cprime{$'$} \def\cprime{$'$} \def\cprime{$'$}
  \def\cprime{$'$} \def\cprime{$'$} \def\cprime{$'$} \def\cprime{$'$}
  \def\cprime{$'$} \def\cprime{$'$} \def\cprime{$'$} \def\cprime{$'$}
  \def\cprime{$'$} \def\cprime{$'$} \def\cprime{$'$} \def\cprime{$'$}
  \def\cprime{$'$} \def\cprime{$'$} \def\cprime{$'$} \def\cprime{$'$}
  \def\cprime{$'$} \def\cprime{$'$} \def\cprime{$'$} \def\cprime{$'$}
  \def\cprime{$'$} \def\cprime{$'$} \def\cprime{$'$} \def\cprime{$'$}
  \def\cprime{$'$} \def\cprime{$'$}
\providecommand{\bysame}{\leavevmode\hbox to3em{\hrulefill}\thinspace}
\providecommand{\MR}{\relax\ifhmode\unskip\space\fi MR }
\providecommand{\MRhref}[2]{%
  \href{http://www.ams.org/mathscinet-getitem?mr=#1}{#2}
}
\providecommand{\href}[2]{#2}

\end{document}